\documentclass[fleqn,11pt,twoside]{article}

\usepackage{amsthm,amsthm,amssymb, color, xcolor,epsfig, graphics, subfigure}

\usepackage{amsmath, graphicx, latexsym, lscape}

\makeatletter
\newcommand{\copyrightnote}[2]{{\renewcommand{\thefootnote}{}
 \footnotetext{\small\it
\begin{flushleft}
 \copyright \ #1   #2  
\end{flushleft}}}}

\newcommand{\Name}[1]{\begin{flushleft}
                       \LARGE \bf #1
                       \end{flushleft}\vspace{-3mm}}

\newcommand{\Author}[1]{\begin{flushleft}
                       \it #1 \end{flushleft}}

\newcommand{\Address}[1]{\begin{flushleft}
                       \it #1 \end{flushleft}}

\newcommand{\Date}[1]{\begin{flushleft}
                      \small  \it #1 \end{flushleft}}

\newcommand{\evenhead}{Author \ name}
\newcommand{\oddhead}{Article \ name}

\renewcommand{\@evenhead}{
\hspace*{-3pt}\raisebox{-15pt}[\headheight][0pt]{\vbox{\hbox to \textwidth
{\thepage \hfil \evenhead}\vskip4pt \hrule}}}
\renewcommand{\@oddhead}{
\hspace*{-3pt}\raisebox{-15pt}[\headheight][0pt]{\vbox{\hbox to \textwidth
{\oddhead \hfil \thepage}\vskip4pt\hrule}}}
\renewcommand{\@evenfoot}{}
\renewcommand{\@oddfoot}{}

\long\def\@makecaption#1#2{%
  \vskip\abovecaptionskip
  \sbox\@tempboxa{\small \textbf{#1.}\ \ #2}%
  \ifdim \wd\@tempboxa >\hsize
    {\small \textbf{#1.}\ \ #2}\par
  \else
    \global \@minipagefalse
    \hb@xt@\hsize{\hfil\box\@tempboxa\hfil}%
  \fi
  \vskip\belowcaptionskip}

\newcommand{\JNMPnumberwithin}[3][\arabic]{%
  \@ifundefined{c@#2}{\@nocounterr{#2}}{%
    \@ifundefined{c@#3}{\@nocnterr{#3}}{%
      \@addtoreset{#2}{#3}%
      \@xp\xdef\csname the#2\endcsname{%
        \@xp\@nx\csname the#3\endcsname .\@nx#1{#2}}}}%
}

\newcommand{\resetfootnoterule} {
  \renewcommand\footnoterule{%
  \kern-3\p@
  \hrule\@width.4\columnwidth
  \kern2.6\p@}
}

\renewcommand{\footnoterule}{}

\makeatother

\theoremstyle{definition}

\begin{document}

\renewcommand{\evenhead}{ {\LARGE\textcolor{blue!10!black!40!green}{{\sf \ \ \ ]ocnmp[}}}\strut\hfill C M Viallet}
\renewcommand{\oddhead}{ {\LARGE\textcolor{blue!10!black!40!green}{{\sf ]ocnmp[}}}\ \ \ \ \  Calculation of the algebraic entropy of a map}

\thispagestyle{empty}
\newcommand{\FistPageHead}[3]{
\begin{flushleft}
\raisebox{8mm}[0pt][0pt]
{\footnotesize \sf
\parbox{150mm}{{Open Communications in Nonlinear Mathematical Physics}\ \ \ \ {\LARGE\textcolor{blue!10!black!40!green}{]ocnmp[}}
\quad Special Issue 1, 2024\ \  pp
#2\hfill {\sc #3}}}\vspace{-13mm}
\end{flushleft}}

\FistPageHead{1}{\pageref{firstpage}--\pageref{lastpage}}{ \ \ }

\strut\hfill

\strut\hfill

\copyrightnote{The author(s). Distributed under a Creative Commons Attribution 4.0 International License}

\begin{center}
{  {\bf This article is part of an OCNMP Special Issue\\ 
\smallskip
in Memory of Professor Decio Levi}}
\end{center}

\smallskip

\Name{An exercise in experimental mathematics:\\
   calculation of the algebraic entropy of a map}

\Author{Claude M. Viallet}

\Address{ Centre National de la Recherche Scientifique \& Sorbonne
  Universit\'e\\ UMR 7589, 4 Place Jussieu F-75252 Paris CEDEX 05,
  France}

\Date{Received August 15, 2023; Accepted January 2, 2024}

\setcounter{equation}{0}

\begin{abstract}

\noindent 
We illustrate the use of the notion of derived recurrences introduced
earlier to evaluate the algebraic entropy of self-maps of projective
spaces. We in particular give an example, where a complete proof is
still awaited, but where different approaches are in such perfect
agreement that we can trust we get to an exact result. This is an
instructive example of experimental mathematics.
\end{abstract}

\label{firstpage}


\section{Introduction: algebraic entropy of maps}
We deal with birational self-maps of $N$-dimensional projective space
$\mathcal{P}_N$. Such maps are given as polynomial maps of degree $d$
when written in terms of the $N+1$ homogeneous
coordinates. Birationality means that the inverse maps are also
polynomial (not necessarily of the same degree). The iterates can be
evaluated polynomially, and the degree $d_n$ of the $n$th iterate is
uniquely defined once all common factors to the homogeneous
coordinates are removed.

The algebraic entropy \cite{BeVi99,Vi98} is defined from the sequence of
degrees $\{d_n\}$ by
\begin{eqnarray} \label{ae}
  \epsilon = \lim_{n\to \infty} \frac{1}{n} \log(d_n)
\end{eqnarray}

Some authors also use the terminology 'dynamic degree'~\cite{RuSh97}
or 'dynamical degree' (see~\cite{Si07} and references therein)
$\delta= \lim_{n\to \infty} d_n^{1/n}$, of which $\epsilon$ is the
logarithm. Our definition is not without relation to the spectral
radius of the map induced in homology, which already appeared in
previous definitions of entropy~\cite{Yo87,Fr91,Gr03} and to the
notion of complexity introduced in~\cite{Ar90}.

The limit~(\ref{ae}) always exists and is invariant by any birational
change of coordinates. It is an excellent -if not the best- detector
of integrability of birational discrete-time systems: integrability =
vanishing of the entropy.

Beyond this use as an integrability criterion, an important question
is to determine the set of non zero values the entropy can take. On
aspect of the problem is then, given a map, {\em to find explicitly the
 exact value of $\epsilon$}.

This can be achieved - especially in the two dimensional case which is
of paramount importance since it encompasses the discrete Painlev\'e
equations~\cite{Sa01,KaNoYa17} - by a singularity analysis of the
maps: in that case one may eventually define a rational variety
obtained by blow-ups of points, where the maps are
diffeomorphisms, and the induced action on the Picard group gives the
answer (see~\cite{DiFa01,Mc07,Du10}).

We are interested in situations where this approach cannot be used. In
particular, going beyond the ``easy'' two-dimensional case, makes the
singularity analysis much more intricate, see for
example~\cite{BeKi04,CaTa19,GrGu23,Vi19}.

The most elementary thing one can do is to evaluate as many terms as
possible of the sequence degrees, using one's favourite formal
calculus software.

An approximate value of $\epsilon$ may then possibly be obtained by
calculating the successive ratios $d_{n+1}/d_n$. Let us call this low
brow analysis 'Method 0'. It gives an idea of the value of $\epsilon$.

A better approach is to look for a generating function for the
sequence of degrees. This function is by definition
  \begin{eqnarray}
    g(s) = \sum_{k=0}^{k=\infty} d_k\ s^k
    \end{eqnarray}
  This method  ('Method 1') - which at first may look too heuristic - works
  extremely well, {\em and this is due to the underlying algebraic
    structure of the problem}. 
  
  A third method ('Method 2'), introduced in~\cite{Vi14} and expanded
in~\cite{Vi15}, is an analysis of the form of the
  iterates. If the factor structure of the iterates happens to
  stabilise, we may rewrite the maps in a different way, giving
  immediately the value of the entropy.
  
  The main point is actually that the asymptotic behaviour which
  $\epsilon$ measures can be extracted from a finite piece of the
  sequence $\{d_n\}$. The fundamental reason is the fact that {\em {
      more than often, the sequence $\{d_n\}$ verifies a finite
      recurrence relation}}.  In addition this recurrence relation has
  integer coefficients, yielding for the entropy the remarkable
  property that it is the logarithm of an algebraic integer. This last
  property was conjectured~\cite{BeVi99} to be true for all birational
  self-maps of projective spaces, and its generality is now
  questioned\cite{BeDiJoKr21}.

  The possible drop of degrees of the iterates (meaning that $d_n$ is
  strictly lower than $d_1^n$), the nature of the generating function
  of the sequence of degrees, as well as the stabilisation phenomenon
  of the form of the iterates are {\em all footprints of the
    singularity structure of the iterations}. Our point is that -
  apart from the simplest two dimensional case, as is the warm-up
  example we have chosen, a full analysis of this structure is often
  difficult for the higher dimensional maps. We stick to an
  experimental approach, using our calculation tools to their limits,
  and do not embark into the singularity analysis, leaving it to
  further study..

  \section{Warm-up: a well studied two dimensional map}

  As a good example of what can be done  for self-maps of $\mathcal{P}_2$, we start from the prototype of algebraically stable (aka ``confining'') map with positive entropy given in~\cite{HiVi98}.
  \begin{eqnarray}
    \label{phihivi}
\varphi: [ x,y,z] \longrightarrow [{x}^{3}+a{z}^{3}-y{x}^{2},{x}^{3},{x}^{2}z]
\end{eqnarray}
which is the transcription as a map in $\mathcal{P}_2$ of  the simple order 2 recurrence
\begin{eqnarray}
  u_{n+1} + u_{n-1} = u_n +\frac{a}{u_n^2}
  \label{hvorig}
\end{eqnarray}

  This map has been shown to have positive entropy by various methods,
among which the construction of a rational surface over $\mathcal{P}_2$ where
the singularities are resolved~\cite{Ta01,Ta01b}. The lift of the map to
the Picard group of this variety is a linear map whose maximal
eigenvalue  (spectral radius)  gives the entropy.

A possible first step is to  calculate the beginning of the sequence of degrees:
\begin{eqnarray}
  \{d_n\} =  1,3,9,27,73,195,513,1347,3529,9243,24201,63363, \dots
  \label{hvliste}
  \end{eqnarray}
and then extract as much as information from this limited amount of data.

Method 0: The most naive - but already useful - thing to do is to see how
the ratio $d_{n+1}/d_n$ evolves with $n$:
\begin{center}
\begin{tabular}{|c|c|}
      \hline
      n & $d_{n+1}/d_n$ \\
      \hline
      1& 3. \\ \hline
      2 &3. \\ \hline
      3 & 3. \\ \hline
      4 & 2.703703... \\ \hline
      5 & 2.671232... \\ \hline
      6 & 2.630769... \\ \hline
      7 & 2.625730... \\ \hline
      8 & 2.619896... \\ \hline
      9 & 2.619155... \\ \hline
      10 & 2.618305... \\ \hline
      11  &  2.618197... \\ \hline
      \end{tabular}
\end{center}
Clearly these numbers point to a value of the order of $\log(2.618...)$
for the entropy. This is not enough if we want the exact value, but it
already tells us that the entropy is not vanishing. In other words the
recurrence does not fall into the integrable class~\cite{HiVi98}.

Method 1: It is possible to fit the sequence (\ref{hvliste}) with the {\em{rational generating function}}
    \begin{eqnarray}
g(s) = \frac{3 \textit{s}^{3}+1}{\left(1-\textit{s}\right) \left(1+\textit{s} \right) \left(\textit{s}^{2}-3 \textit{s} +1\right)}
      \end{eqnarray}
The entropy is the logarithm of the  inverse of the smallest  modulus of the poles of $g$, since this is what governs the growth of the Taylor expansion of $g(s)$.
This gives  $\epsilon=
\log((3+\sqrt{5})/2) \simeq \log(2.618033988...) $

To go beyond the approximation of Method 0 and the heuristic nature of
{Method 1}, go to Method 2, keeping in mind that it may imply much heavier
calculations.

For the map~(\ref{phihivi}) the form of the iterates does stabilise to
the pattern
\begin{eqnarray}
p_k=[  A_{k-3}^3\, A_k,  A_{k-4}\, A_{k-1}^3, z\, A_{k-3}^2\, A_{k-2}^2\, A_{k-1}^2 ]
\label{stablehv}
\end{eqnarray}
and the recurrence relation between the blocks $A_k$ is just
\begin{eqnarray}
\label{hirohv}
A_k^3\, A_{k-3}^3+a\, z^3\, A_{k-1}^6\, A_{k-2}^6 - A_{k-1}^3\,
A_{k-4}\, A_{k}^2 =\, A_{k-3}^2\, A_{k-2}^3\, {\bf A_{k+1}}
\end{eqnarray}
We call~(\ref{hirohv}) the {\em{derived recurrence}} of the original
one~(\ref{hvorig}).  Notice that it is the same as equ (4.6)
of~\cite{Ho07}, where it was obtained by a different
approach. Recurrence~(\ref{hirohv}) extends over a string of length
$6$.  This relation is not quadratic nor multi-linear, but it allows
to prove (\ref{stablehv}), providing the recurrence condition on the
degrees of the iterates of $\varphi$, and the value of the entropy
$\epsilon=\log( (3+\sqrt{5})/2)$ (same as above).

The properly conducted singularity analysis, blowing up enough points
of $\mathcal{P}_2$, and looking at the induced map on the Picard group of the
rational variety constructed in this way, confirms the value
$\epsilon= \log( (3+\sqrt{5})/2)$\cite{Ta01,Ta01b}.

So all Methods 0,1,2, and the complete singularity analysis agree
perfectly. Method 0 is approximate but useful, Method 1 is providing
us with an educated guess, and suggests a candidate for the value of
$\epsilon$.  This value turns out to be exact, as one can prove using
Method 2, or with the full desingularisation of the map.

Remark: The sequence (\ref{hvliste}) is registered in the On-Line
Encyclopedia of Integer Sequences (https://oeis.org/) under number A084707.

At this point it is interesting to notice that the recurrence
relation~(\ref{hirohv}) belongs the 'Somos-like' family
(see~\cite{So89,EkZe14,Ze} and
https://faculty.uml.edu//jpropp/somos/history.txt ). Although
$A_{k+1}$ is given as a rational fraction in terms of the previous
$A$'s, the recurrence has the so called Laurent
property~\cite{Vi15,KaMaTo18}. Moreover, if one launches the
recurrence with appropriate polynomial initial conditions, {the
  values one obtains are, {\em by construction}, multivariate
  polynomials}.
  
\section{ A more challenging map}

Monomial maps are known to behave in a particular way, as far as the
sequence of degrees of their iterates is concerned. Birational
monomial maps have an entropy which is the logarithm of an algebraic
integer, but the generating function of the sequence is not
necessarily rational~\cite{HaPr05}.

We will now examine a map acting in dimension larger than two.
Consider the  recurrence  of order 4~\cite{Gu22}:
\begin{eqnarray}
  x_{n+1} = \frac{x_n x_{n-2}}{x_{n-1} ( 1- x_{n-1}) x_{n-3}} 
  \end{eqnarray}

This recurrence defines  a  birational (almost monomial) map in $\mathcal{P}_4$.
\begin{eqnarray}
  \label{f-map}
   \varphi:   [x,y,z,u,t] \rightarrow  [x\,z\,t^2, x\,u\,y\,(t-y), u\,y^2\,(t-y), u\,y\,(t-y)\,z, u\,y\,(t-y)\,t]
\end{eqnarray}
with inverse
\begin{eqnarray}
\psi:  [x,y,z,u,t] \rightarrow [ x\,y\, z \left(t-z \right),  x\, z^2  \left(t-z \right),x\,z\,u \left(t-z \right),  y\, u\, t^2, x\, z  \left(t -z \right) t]
\end{eqnarray}

The direct calculation of the sequence of degrees yields
\begin{eqnarray}
  \label{degrees_fmap}
\{d_n\}=&&1, 4, 5, 9, 11, 16, 21, 30, 43, 61, 86, 120, 168, 234, 329, 459, 645, 902, 
1267,
1771, \nonumber\\
&& 2484, 3476, 4871, 6822, 9555, 13384, 18745, 26256, 36774, 51507, 72143,
101043, \nonumber\\
&& 141524, 198223, 277633, 388864, 544644, 762846, 1068451,
1496494, 
2096019, \nonumber \\
&& 2935716, 4111826, 5759091, 8066291, 11297797, 
15823888,
22163239, 31042218, \nonumber \\
&& 43478302, 60896502, 85292724, 119462566, 1677321393,234353404,328239604,
 \nonumber \\ &&  \dots
\end{eqnarray}

\subsection{Method 0}
We have  the following sequence of ratios $d_{n+1}/d_n$:
\begin{eqnarray*}
 4.,  &&  1.250000000,    1.800000000,    1.222222222,    1.454545455,    1.312500000,  1.428571429,  \\ &&    1.433333333,    1.418604651,    1.409836066,    1.395348837,    1.400000000,  1.392857143,   \\ &&   1.405982906,    1.395136778,    1.405228758,    1.398449612,    1.404656319,  1.397790055,   \\ &&   1.402597403,    1.399355878,    1.401323360,    1.400533771,    1.400615655,   1.400732601, \\ &&    1.400552899,    1.400693518,    1.400594150,    1.400636319,    1.400644573, 1.400593266,  \\ &&     1.400631414,    1.400631695,    1.400609415,    1.400640414,    1.400602781,  1.400632340,   \\ &&   1.400611657,    1.400620150,    1.400619715,    1.400615166,    1.400621177,  1.400616417,  \\ &&    1.400618778,    1.400618574,    1.400617129,    1.400619051,    1.400617392,   1.400618409,  \\ &&   1.400618221,    1.400617789,    1.400618487,    1.400617772,    1.400618294,   1.400618034,  \\ &&  
\dots  
  \end{eqnarray*}
There is no doubt that the sequence converges, and that the numerical
value of the entropy is around $\log(1.400618...)$, but we want to
have its exact value.

\subsection{Method 1}
It is not possible to fit the sequence~(\ref{degrees_fmap}) with a
satisfactory rational generating function. This could come from two
different reasons: either the generating function is not rational,
either the information we have is not sufficient. Unfortunately, the
explicit calculation of the degrees cannot go much further due to the
practical limitations of the formal calculus software we have at
hand...

\subsection{Method 2}
The images $p_n$ of the generic starting point $p_0=[x,y,z,u,t]$ are
made of products of factors $B_k$ with $B_1=t-y$ and $B_2=t-x$, and
the further $B_k$'s are the proper transforms\cite{Sh77} of $B_{k-1}$.

The first coordinate of $p_n$ is a product of some $B_j$'s with
various powers,  $B_n$ not appearing, and some adventive monomials
in $ x,y,z,u,t$, which remain of low degree.

The other four coordinates are of the same form, but all contain
$B_n$ with power 1.

The outcome\footnote{Maple calculation up to order 17, insufficient,
and then calculation using the V. Shoup's NTL C++ library~\cite{NTL}}
is that the form indeed stabilises after order 27,
and remains unchanged up to the maximum order we were able to
reach. To give an idea the 44th iterate looks like
\begin{eqnarray*}
\hskip -1truecm
&&p_{44}:= [y\, u^3\, t^3\, B_{18}\, B_{19}\, B_{20}\, B_{23}\, B_{24}^2\, B_{25}^2\, B_{26}^3\, B_{27}^3\, B_28\, B_{31}^2\, B_{32}\, B_{33}^3\, B_{34}^4\, B_{35}^2\, B_{38}\, B_{39}\, B_{40}^2\, B_{41}^3\, B_{42}^2,\\ &&
x^2\, z\, u^4\, t^2\, B_{18}\, B_{19}\, B_{22}^2\, B_{23}^2\, B_{24}\, B_{25}^2\, B_{26}^3\, B_{27}\, B_{29}^2\, B_{30}^2\, B_{31}\, B_{32}\, B_{33}^4\, B_{34}^2\, B_{36}\, B_{37}^2\, B_{38}\, B_{39}\, B_{40}^3\, B_{41}^2\, B_{44},\\ &&
x^2\, y\, u^2\, t^3\, B_{18}\, B_{21}^2\, B_{22}^3\, B_{23}\, B_{24}\, B_{25}^2\, B_{26}\, B_{28}^2\, B_{29}^4\, B_{30}\, B_{31}\, B_{32}^2\, B_{33}^2\, B_{35}\, B_{36}^3\, B_{37}^2\, B_{38}\, B_{39}^2\, B_{40}^2\, B_{43}\, B_{44},\\&&
x^2\, y^2\, t^3\, B_{20}^2\, B_{21}^3\, B_{22}^2\, B_{23}\, B_{24}\, B_{27}^2\, B_{28}^4\, B_{29}^3\, B_{30}\, B_{31}^2\, B_{34}\, B_{35}^3\, B_{36}^3\, B_{37}^2\, B_{38}^2\, B_{39}\, B_{42}\, B_{43}\, B_{44},\\ &&
x\, y\, u^2\, B_{18}\, B_{19}\, B_{20}\, B_{21}\, B_{22}\, B_{25}\, B_{26}^2\, B_{27}^2\, B_{28}^2\, B_{29}^2\, B_{33}^2\, B_{34}^2\, B_{35}^2\, B_{36}^2\, B_{37}\, B_{40}\, B_{41}\, B_{42}\, B_{43}\, B_{44}];
\end{eqnarray*}

and the equation giving $B_{44}$ is
\begin{eqnarray}
  \label{equ44}
equ44: &\{& B_{17}\,{\bf B_{44}} =   
  B_{19}\,B_{20}\,B_{26}\,B_{27}^2\,B_{34}^2\,B_{35}\,B_{41}\,B_{42}
\nonumber  \\ &-&
  B_{21}\,B_{22}^2\,B_{23}\,B_{24}\,B_{25}\,B_{29}^2\,B_{30}\,B_{31}\,B_{32}^2\,B_{36}\,B_{37}\,B_{38}\,B_{39}^2\,B_{40}\,t^3\,x\}
\end{eqnarray}

This is a recurrence of order 27, and it gives for the rate of growth of the degrees  the largest root of
\begin{eqnarray}
  r^{44}-r^{42}-r^{41}-r^{35}-2 r^{34}-2 r^{27}-r^{26}-r^{20}-r^{19}+r^{17}
  \label{rateB}
  \end{eqnarray}
which happens to be  1.400618098... to be compared with the approximate value  $\epsilon \simeq \log(1.400618...)$ we had from the original sequence of degrees.

Remark: In (\ref{rateB}) we neglected the factors  in $x,y,z,u,t$ appearing in (\ref{equ44}). These factors are present and will be taken into account in the next section. The main point is that {\em they do not affect the value of the entropy}.

Happily enough we have a candidate for the exact value of the
entropy. It is the logarithm of an algebraic integer.

Notice that the derived recurrence if of much larger order than the
initial one.

\section{Why did Method 1 not work?}
The monomial factors appearing in $p_n$ have some structure:

Setting
\begin{eqnarray*}
  f(k)= 
  {\frac
    {
      B_{k-21}\, B_{k-20}^{2}  B_{k-19}\,  B_{k-18}\,  B_{k-17}\, B_{k-13}^{2}  B_{k-12}\,  B_{k-11}\, B_{k-10}^{2}  B_{k-6}\,  B_{k-5}\,  B_{k-4}\, B_{k-3}^{2}  B_{k-2}
    }
    {
      B_{k-23}\,  B_{k-22}\,   B_{k-16}\, B_{k-15}^{2} B_{k-8}^{2}  B_{k-7}\,  B_{k-1}\,{ B_k}
    }
  },
\end{eqnarray*}
the $k$-th iterate $p_k$  of the generic point $p_0=[x,y,z,u,t]$ then reads
\begin{eqnarray}
  \label{pk}
  p_k \simeq [ \alpha_k \cdot f(k), \beta_k\cdot f(k-1), \gamma_k \cdot f(k-2), \delta_k\cdot f(k-3), \eta_k]
\end{eqnarray}
where $\simeq$ means equality up to a common factor.

The adventive factors $\alpha_k, \beta_k, \gamma_k, \delta_k, \eta_k$
are simple monomials in $x,y,z,u,t$.  Two periods appear for these
factors, namely 7 and 32. Defining
\begin{eqnarray}
\rho_k &=& x^{X_{\rho}(k)} \,  y^{Y_{\rho}(k)} \,  z^{Z_{\rho}(k)} \,  u^{U_{\rho}(k)} \,
t^{T_{\rho}(k)} , \quad   \rho=\alpha,\beta,\gamma,\delta,\eta.
\end{eqnarray}

The powers $X_\rho$, $U_\rho$ and $T_\rho$ are periodic with period
32, and $Y_\rho$ and $Z_\rho$ are periodic with period 7. There exist in addition
simple relations between $Y_\rho$ and $Z_\rho$ as well as between  $U_\rho$ and
$X_\rho$, $ \rho=\alpha,\beta,\gamma,\delta,\eta$.
\begin{eqnarray}
  Y_{\rho}(k) = Z_{\rho}(k+2) \quad \mbox {and} \quad  U_{\rho}(k) = X_{\rho}(k+4), \quad   \rho=\alpha,\beta,\gamma,\delta,\eta.
  \end{eqnarray}
We give here the  values of $X$, $Y$, and $T$  on one period starting at $k=1$
\begin{eqnarray*}
X_\alpha&=&      [1,0,0,0,0,0,2,3,2,1,1,0,0,2,4,3,1,2,0,0,1,3,3,2,2,1,0,0,1,1,1,1],\\
X_\beta& =&      [1,1,0,0,1,0,0,2,3,1,1,2,1,0,2,4,1,1,2,2,0,1,3,2,1,2,2,0,0,1,1,0],\\
X_\gamma&=&   [0,1,1,0,1,1,0,0,2,2,1,2,3,1,0,2,2,1,1,4,2,0,1,2,1,1,3,2,0,0,1,0],\\
X_\delta&=&      [0,0,1,1,1,1,1,0,0,1,2,2,3,3,1,0,0,2,1,3,4,2,0,0,1,1,2,3,2,0,0,0],\\
X_\eta&=&        [0,0,0,0,1,1,1,1,1,0,0,1,2,2,2,2,0,0,0,2,2,2,2,1,0,0,1,1,1,1,1,0],\\
Y_\alpha&=&   [0, 1, 2, 1, 0, 1, 0], \qquad  
Y_\beta =  [1, 0, 1, 2, 0, 0, 1],\qquad
Y_\gamma=[2, 1, 0, 1, 1, 0, 0],\\  
Y_\delta &=&   [1, 2, 1, 0, 0, 1, 0],\qquad
Y_\eta =     [1, 1, 1, 1, 0, 0, 0],\\  
T_\alpha&=&     [2,4,3,1,1,0,0,2,6,6,3,3,1,0,0,4,6,5,3,3,0,0,2,4,4,3,3,1,0,0,0,0],\\
T_\beta& =&     [0,2,4,2,1,2,2,0,2,6,4,2,3,4,0,0,4,6,2,3,4,2,0,2,4,2,2,3,2,0,0,0],\\
T_\gamma&=&  [0,0,2,3,2,2,4,2,0,2,4,3,2,6,4,0,0,4,3,2,4,6,2,0,2,2,1,2,4,2,0,0],\\
T_\delta&=&     [0,0,0,1,3,3,4,4,2,0,0,3,3,5,6,4,0,0,1,3,3,6,6,2,0,0,1,1,3,4,2,0],\\
T_\eta&=&       [1,1,1,0,0,1,3,3,3,3,1,0,0,3,3,3,3,3,0,0,1,3,3,3,3,1,0,0,1,1,1,1].
\end{eqnarray*}

All this means that there is a global period of $224  = 7 \times  32$ for these factors.

The $B_k$ verify the recurrence relation of order 27
\hskip -2truecm
\begin{eqnarray} \label{recurB}
\hskip -2truecm
&&
  \mu_k \cdot B_{k-23}\, B_{k-22}^2 B_{k-21}\, B_{k-20}\, B_{k-19}\, B_{k-15}^2 B_{k-14}\, B_{k-13}\, B_{k-12}^2 B_{k-8}\, B_{k-7}\,   B_{k-6}\, B_{k-5}^2 B_{k-4} \nonumber  \\
  &&   + \; \nu_k \cdot B_{k-25}\, B_{k-24}\, B_{k-18}\, B_{k-17}^2 B_{k-10}^2 B_{k-9}\, B_{k-3}\, B_{k-2}
-B_{k-27}\, {\bf B_{k}}=0
  \end{eqnarray}
where $\mu_k$ and $\nu_k$ are monomials in $x,y,z,u,t$ with period
224. They can be evaluated from the explicit expression of $p_k$.

The generating function of both sequences ( $g_B$ for the degrees of the $B$'s, and $g_p$ for the $p_k$'s)  are then readily obtained.

\begin{eqnarray*}
  g_B= {\frac { P}
{ \left( s-1  \right)  \left( s+1 \right)  \left( {s}^{6}+{s}^{5}+{s}^{4}+{s}^{3}+{
s}^{2}+s+1 \right)  \left( {s}^{4}+1 \right)  \left( {s}^{8}+1
 \right)  \left( {s}^{16}+1 \right) Q }},
\end{eqnarray*}
with
\begin{eqnarray*}
 P & = &  s^{58}+s^{57}-s^{56}-2 s^{55}-2 s^{52}-2 s^{51}+s^{50}+s^{49}-2 s^{48}-2 s^{47}+s^{46}-3 s^{44}-2 s^{43}\\ && +2 s^{42}+s^{41}-3 s^{40}-2 s^{39}+2 s^{38}-4 s^{36}-2 s^{35}+3 s^{34}+s^{33}-3 s^{32}-s^{31}+2 s^{30} \\ && -3 s^{29}-7 s^{28}-3 s^{27}+s^{26}-3 s^{25}-5 s^{24}-s^{23}+s^{22}-4 s^{21}-6 s^{20}-2 s^{19}-4 s^{17} \\ && -4 s^{16}-s^{15}-4 s^{13}-5 s^{12}-2 s^{11}-s^{10}-4 s^{9}-3 s^{8}-3 s^{5}-3 s^{4}-s^{3}-s^{2}-2 s -1,
\end{eqnarray*}
and
\begin{eqnarray*}
 Q =  {s}^{22}-{s}^{20}-{s}^{19}
-{s}^{17}+{s}^{15}+{s}^{14}-{s}^{13}-{s}^{12}-{s}^{10}-{s}^{9}+{s}^{8}
+{s}^{7}-{s}^{5}-{s}^{3}-{s}^{2}+1.
\end{eqnarray*}
We also have
\begin{eqnarray*}
  \label {gp}
  g_p ={\frac {R} { \left(1- s \right) \left( s+1 \right) \left(
      {s}^{2}+1 \right) \left( {s}^{4} +1 \right) \left(
      {s}^{6}+{s}^{5}+{s}^{4}+{s}^{3}+{s}^{2}+s+1 \right) \left(
      {s}^{8}+1 \right) \left( {s}^{16}+1 \right) Q }},
\end{eqnarray*}
with
\begin{eqnarray*}
  R& = & 4\,{s}^{59}+9\,{s}^{58}+14\,{s}^{57}+16\,{s}^{56}+18\,{s}^{55
  }+15\,{s}^{54}+13\,{s}^{53}+10\,{s}^{52}+13\,{s}^{51}+15\,{s}^{50}
  \\ &&
    +19 \,{s}^{49}+21\,{s}^{48}+22\,{s}^{47}+19\,{s}^{46}+15\,{s}^{45}+11\,{s}
^{44}+11\,{s}^{43}+15\,{s}^{42}+18\,{s}^{41}\\ && +23 \,{s}^{40}
+23\,{s}^{39}
+22\,{s}^{38}+15\,{s}^{37}+11\,{s}^{36}+6\,{s}^{35}+6\,{s}^{34}+6\,{s}
^{33}+12\,{s}^{32}+15\,{s}^{31}
\\ && 18\,{s}^{30} 
+15\,{s}^{29}
+13\,{s}^{28}
+7\,{s}^{27}+6\,{s}^{26}+5\,{s}^{25}+10\,{s}^{24}+13\,{s}^{23}+20\,{s}
^{22}+21\,{s}^{21} \\ &&  +22\,{s}^{20}
+17\,{s}^{19}
+14\,{s}^{18}
+10\,{s}^{17}
+10\,{s}^{16}+13\,{s}^{15}+17\,{s}^{14}+20\,{s}^{13}+20\,{s}^{12}
+18\,{s}^{11} \\ && +14\,{s}^{10}
+12\,{s}^{9}+9\,{s}^{8}
+11\,{s}^{7}
+13\,{s}^{6}+
16\,{s}^{5}+15\,{s}^{4}+13\,{s}^{3}+9\,{s}^{2}+5\,s+1.
  \end{eqnarray*}

The entropy is the log of the inverse of the smallest root of the
polynomial $Q$ (which also happens to be its largest root)
approximately $ \log (1.400618098)$ in perfect agreement with the
numerical evaluation obtained from the explicit calculation of the 55
first terms of the sequence of degrees of the iterates of the map.

The degrees of the numerator and the denominator (respectively 59 and
60) of the generating function $g_p$ indicate that we would have
needed to evaluate the degree of the first 119 iterates to use Method
1.  This is beyond the capabilities of the presently available formal
calculus software.

Notice that the denominator of $g_p$ is just $Q\; (s^7-1)\,(s^{32}-1)/
(s-1)$, and the two periods 7 and 32 appear there naturally (they are
present as well in $g_B$).

All this is not quite a proof, but we can bet we have the exact value
of the entropy, {\em and this value is once more the logarithm of an
  algebraic integer!}.

\section{Conclusion}
The derived recurrence~(\ref{recurB}) has two remarkable properties,
having to do with Laurent/'Somos like' characteristics.

The first one is obtained by construction: if we take as initial
conditions the factors $B_k$ obtained from the iterates of the generic
point $[x,y,z,u,t]$ of $\mathcal{P}_4$, {\em all further $B$'s are
  multivariate polynomials}.

Moreover -and this is another {\em experimental} fact- if we start
from $B_i=1, i=1..27$ then again all further $B's$ are polynomials in
$(x,y,z,u,t)$, and of course integers if $(x,y,z,u,t)$ are themselves.
This comes from the fact that recurrence~(\ref{recurB}) {\em has the
  Laurent property for arbitrary $[x,y,z,u,t]$ }, as is easy to check
explicitly on the first iterates. The proof will come later.

In summary, the 'derivation' process of recurrences provides us with a
{\em factory of Somos like recurrences}, keeping in mind that it is
not a mere change of coordinates. It is a complete change of description.

This raises a number of questions:
\begin{itemize}
  \item{ The process defines sequences of multivariate
    polynomials. What are the properties of these polynomials?}
\item{ Is it possible to predict the order of the derived recurrence?
  We gave two examples. In the first one the original recurrence was
  of order 2 and the derived one of order 5, and for the latter the
  original order was 4 and the new one is 27.}
\item{ Iterating the derivation process could produce more recurrences
  but may as well reach a fixed point. This should be investigated. In
  fact, when applied to the Somos-4 recurrence, it just reproduces
  Somos-4 with a periodic decoration similar to the factors $\mu_k$
  and $\nu_k$ seen in~(\ref{recurB})}\footnote{I would like to thank
M. Somos for drawing my attention to this case, and fruitful further
exchanges, as well as for  signalling a sign mistake in the first draft of
this paper.}.
\item{ The stabilisation of the form of the iterates is necessary for
  the existence of a derived recurrence. How can we characterise the
  systems having this property? This question is crucial since it
  ensures that the entropy is the log of an algebraic integer.}
  \item{ In the specific case of discrete integrable systems
    (vanishing entropy) the derived recurrence may take the form of
    the Hirota quadratic relation between $\tau$ functions. The link
    has to be clarified.}
\end{itemize}

There would be much more to say about last item of the previous list:
Discrete Integrable Systems. Among the contributors to the subject
Decio Levi was there from the first day. He was at the origin
(together with Pavel Winternitz, and Luc Vinet) of the series of SIDE
conferences~\cite{SIDE}, a very important series of meeting in the
field and he took an active part, as early as 1994, organising more
than one of the meetings. I have had a long interaction with Decio,
not only at the occasion of these conferences, but also of visits to
Roma Tre (where I made the acquaintance of his then Ph.D. student
Giorgio Gubbiotti, now a collaborator). Decio was always very
supportive, and discussions with him very constructive.  He will be
sorely missed.


\label{lastpage}
\end{document}